\documentclass{article}
\usepackage{latexsym,amsmath,amssymb,amscd}
\usepackage[mathscr]{eucal}
\usepackage[dvips]{color}
\begin{document}

\title{A Remark on the Matter-Vacuum Matching Problem for Axisymmetric Metrics Governed by the Einstein-Euler Equations }
\author{Tetu Makino \footnote{Professor Emeritus at Yamaguchi University, Japan.  e-mail: makino@yamaguchi-u.ac.jp}}
\date{\today}
\maketitle

\newtheorem{Lemma}{Lemma}
\newtheorem{Proposition}{Proposition}
\newtheorem{Theorem}{Theorem}
\newtheorem{Definition}{Definition}
\newtheorem{Remark}{Remark}
\newtheorem{Corollary}{Corollary}
\newtheorem{Notation}{Notation}
\newtheorem{Assumption}{Assumption}

\numberwithin{equation}{section}

\begin{abstract}
Axially symmetric stationary metrics governed by the Einstein-Euler equations for slowly rotating perfect fluids have been constructed in an arbitrarily large bounded domain containing the support of the mass density. However the problem of global prolongation of the metric is still open. On the other hand the so called matter-vacuum matching problem,
particularly as the source problem for the Kerr metric, has been discussed by several authors. This can be regarded as the approach to the same open problem in the opposite direction. We give a remark on this open problem.
\end{abstract}

\section{Introduction}

We consider the Einstein-Euler equations
\begin{equation}
R_{\mu\nu}-\frac{1}{2}g_{\mu\nu}R=\frac{8\pi\mathsf{G}}{\mathsf{c}^4}T_{\mu\nu}
\end{equation}
for the metric
\begin{equation}
ds^2=g_{\mu\nu}dx^{\mu}dx^{\nu}
\end{equation}
with the energy-momentum tensor for the perfect fluid
\begin{equation}
T^{\mu\nu}=(\mathsf{c}^2\rho+P)U^{\mu}U^{\nu}-Pg^{\mu\nu}.
\end{equation}
Here the Greek letters $\mu, \nu$ stand for $0,1,2,3$. $\mathsf{G}$ is the gravitational constant and $\mathsf{c}$ is the speed of light. As for the meanings of other symbols, see
\cite{LandauL}. We assume the pressure $P$ is a given function of the density $\rho$ and suppose \\

{\bf (A) : $P$ is a given smooth function of $\rho >0$ such that $0<P, 0<dP/d\rho <\mathsf{c}^2$ for $\rho >0$, and there is a smooth function $\Lambda$ which is analytic at $0$ such that $\Lambda(0)=0$ and
\begin{equation}
P=\mathsf{A}\rho^{\gamma}(1+\Lambda(\mathsf{A}\rho^{\gamma-1}/\mathsf{c}^2)) \quad\mbox{for}\quad \rho >0.
\end{equation}
Here $\mathsf{A}, \gamma$ are positive constants and $6/5 < \gamma <2$.}\\

In the article \cite{asEE}, we have constructed an axisymmetric metric
\begin{align}
ds^2=&e^{2F'(\varpi, z)}(\mathsf{c}dt+A'(\varpi, z)d\phi')^2-
e^{-2F'(\varpi, z)}\Big[
e^{2K'(\varpi, z)}(d\varpi^2+dz^2)+ \nonumber \\
&\Pi(\varpi, z)^2(d\phi')^2\Big],\qquad\mbox{while}\quad  \phi'=\phi-\Omega t,
\end{align}
with the 4-velocity field
\begin{equation}
U^{\mu}\frac{\partial}{\partial x^{\mu}}=\frac{1}{\mathsf{c}F'(\varpi, z)}\Big(\frac{\partial}{\partial t}+\Omega\frac{\partial}{\partial \phi}\Big),
\end{equation}
and the density distribution
\begin{equation}
\rho=\Big(\frac{\gamma-1}{\mathsf{A}\gamma}\Big)^{\frac{1}{\gamma-1}}
(u\vee 0)^{\frac{1}{\gamma-1}}\Big(1+O\Big(\frac{u}{\mathsf{c}^2}\Big)\Big), \quad u=u_{\mathsf{N}}+O\Big(\frac{u_{\mathsf{O}}^2}{\mathsf{c}^2}\Big)
\end{equation}
in a bounded domain $\mathfrak{D}=\{ (\xi^1,\xi^2, \xi^3) \in \mathbb{R}^3 | r:=\sqrt{\varpi^2+z^2}<R\}$,
where $\xi^1=\varpi\cos\phi, \xi^2=\varpi\sin\phi, \xi^3=z$,  provided that 
$u_{\mathsf{O}}/\mathsf{c}^2$ is sufficiently small. Here 
$R$ is an arbitrarily fixed large number, and $u_{\mathsf{N}}(r,\zeta), \zeta:=z/r, 
u_{\mathsf{O}}=u(O),$
is the solution of the problem of axially symmetric rotating gaseous stars for the non-relativistic Euler-Poisson equations such that $\{u>0\}=\{\rho >0\}=\{r <r_+(\zeta)\}$,
the function
$ \zeta \mapsto r_+(\zeta): [-1,1] \rightarrow ]0, R/2[$ being a continuous function which gives the physical vacuum boundary. 
$F', A', K', \Pi $ and $u$ as functions of $(\xi^1,\xi^2,\xi^3)$ belong to
$C^{2,\alpha}(\overline{\mathfrak{D}}), 0<\alpha <1$. 
We have assumed that $\Omega^2/2\pi\mathsf{G}\rho_{\mathsf{O}}, \displaystyle \rho_{\mathsf{O}}:=\Big(\frac{\gamma-1}{\mathsf{A}\gamma}\Big)^{\frac{1}{\gamma-1}}u_{\mathsf{O}}^{\frac{1}{\gamma-1}},$ is sufficiently small. \\

Now, in this situation, the support of $\rho$ is $\overline{\mathfrak{R}}:=\{ r \leq r_+(\zeta)\}$ is a compact subset of the domain $\mathfrak{D}$, on which the metric has been constructed, and $\mathfrak{D} \setminus \overline{\mathfrak{R}}$ is a vacuum region. However there remains the open problem:\\

{\bf What happens when we prolong the vacuum axisymmetric metric outside $\mathfrak{D}$ ?} \\

{\bf Can the prolongation be an asymptotic flat one as $r \rightarrow +\infty$ ?}\\

For the sake of comparison, let us recall the result on spherically symmetric metric. Actually a spherically symmtric static metric
\begin{equation}
ds^2=e^{2F(r)}\mathsf{c}^2dt^2
-e^{2H(r)}dr^2-r^2(d\vartheta^2+\sin^2\vartheta d\phi^2) \label{ssInt}
\end{equation}
is given by solving the Tolman-Oppenheimer-Volkoff equation. Let 
$\{\rho >0\}=\{r <r_+\}$, where $r_+$ is supposed to be finite, which is the case 
if $4/3 <\gamma$ or if
$u_{\mathsf{O}}/\mathsf{c}^2$ is sufficiently small. Anyway, in this situation, we can find a Schwarzschild metric described by the coordinates $(R^{\sharp}, \vartheta, \phi)$:
\begin{align}
ds^2=&\Big(1-\frac{2\mathsf{G}m_+}{\mathsf{c}^2R^{\sharp}}\Big)\mathsf{c}^2(dt^{\sharp})^2-
\Big(1-\frac{2\mathsf{G}m_+}{\mathsf{c}^2R^{\sharp}}\Big)^{-1}(dR^{\sharp})^2 \nonumber \\
&-(R^{\sharp})^2(d\vartheta^2+\sin^2\vartheta d\phi^2) \label{ssExt}
\end{align}
on the domain $\{ r_+ < R^{\sharp}\}$ such that the coefficients of the patched metric 
\eqref{ssExt} with \eqref{ssInt} at $r=R^{\sharp}=r_+$ are of class $C^2$. See \cite{ssEE}.

Indeed, by Birkhoff's theorem the Schwarzschild metric the only possible spherically symmetric static one in the vacuum region. 

Keeping in mind this spherically symmetric case, and that the Kerr metric is the typical, although not only one, static axially symmetric metric in the vacuum region, we naturally ask: \\

{\bf (Q1): Can we find a Kerr metric described by an appropriate coordinates on the vacuum region
$\{ r_+(\zeta) <r\}$ which is matched to the interior solution constructed on
$\mathfrak{R}=\{\rho >0\}=\{r <r_+(\zeta)\}$ with junction conditions to be expected, say, e.g., with coefficients of class $C^1$?}\\

Let us discuss this matter-vacuum matching problem in this note.

\section{Strategy of W. Roos}

Until now the matter-vacuum matching problem has been discussed under the formulation which is in the opposite direction to the formulation described in the preceding section
as {\bf (Q1)}. Namely, the formulation of the problem discussed by several authors until now is: \\

{\bf (Q2): Given an axisymmetric stationary vacuum metric, in particular, the Kerr metric, on the exterior domain $\mathcal{E}=\mathbb{R}^3\setminus \overline{\mathcal{I}}$ of an interior domain $\mathcal{I}$ with boundary
$\mathcal{S}=\partial\mathcal{E}=\partial\mathcal{I}$, can we find a solution of the field equations such that $\{ \rho >0\}$ or $\{ P >0\}$ coincides with the interior $\mathcal{I}$ on which the metric is constructed to be matched across $\mathcal{S}$ with plausible junction conditions, say, e. g., with coefficients of class $C^1$ to the given exterior  vacuum metric?}\\

Several authors have been discussed this problem of so called `a source of the Kerr field'.
The textbook \cite{PlebanskiK} by J. Pleba\'{n}ski and A. Krasi\'{n}ski, 2006, pp.499-495, says

\begin{quote}
The Kerr solution has been known for more than 40 years now, and from the very beginning
its existence provoked the simple question: what material body could generate such a vacuum field around it? Several authors have tried very hard to find a model of the source, but so far without success. The most promising positive result is that of Roos (1976), who investigated the Einstein equations with a perfect fluid with the boundary condition that the Kerr metric is matched to the solution. All attempts  so far to find an explicit example of a solution failed. The continuing lack of success prompted some authors to spread the suspicion that a perfect fluid source might not exist; rumours about this suspicion were then taken as a serious suggestion. The opinion of one of the present authors (A. K.) is that a bright new idea is needed, as opposed to routine standard tricks tested so far. 
\end{quote}

However the strategy of W. Roos \cite{Roos} seems to be not so promising. In this section we examine this.\\

A rough sketch of the main result of \cite{Roos} giving Roos' strategy is as follows:\\

Let be given a stationary axisymmetric $C^{\infty}$ vacuum filed in a domain $\tilde{\mathcal{E}}$, an analytic equation of state $\rho=\rho(P)$, a time-like hypersurface $\mathcal{S}$ bounding a simply connected domain $\mathcal{I}$, $\mathcal{E}=\mathbb{R}^3\setminus \bar{\mathcal{I}} \subset
\tilde{\mathcal{E}}$, such that on $\mathcal{S}$ matching conditions and some additional conditions hold. Take a coordinate system $x^1,x^2,x^3$ in a neighborhood $\mathcal{V}$ of $\mathcal{S}$ such that
$\{x^1 < 0\}\cap\mathcal{V}=\mathcal{E}\cap\mathcal{V}$,
or, $\{x^1 > 0\}\cap\mathcal{V}=\mathcal{I}\cap\mathcal{V}$.
 Then we can find, in a neighborhood of $\mathcal{S}$, a unique stationary axisymmetric and analytic field for a rotating perfect fluid satisfying the matching conditions on $\mathcal{S}$ by applying the Cauchy-Kovalevskaja theorem on the initial surface $\{x^1=0\}$ so that $\partial P/\partial x^1 >0$ on $x^1>0$. \\

For the details of the matching conditions and additional conditions, see \cite[(4.1)-(4.3), (5.3)]{Roos}. However it is almost everywhere impossible to expect that the Roos' strategy would successfully give the  solution in the whole interior domain $\mathcal{I}$. Let us show this defect. \\

Let us consider the particular case of spherically symmetric solutions. For the sake of simplicity let us take the geometrical unit system in which  $\mathsf{G}=\mathsf{c}=1$.

 Given the Schwarzschild metric 
$$ds^2=\Big(1-\frac{2M}{r}\Big)dt^2-
\Big(1-\frac{2M}{r}\Big)^{-1}dr^2
-r^2(d\vartheta^2+\sin^2\vartheta d\phi^2) $$
on $\tilde{\mathcal{E}}=\{ r >2M\}$ and $\mathcal{S}=\{r=R\}$, where the data $(R, M)$ belongs to the admissible set $\mathcal{A}$ defined by
$$\mathcal{A}:=\{(R,M) | R >0, M>0, 1-\frac{2M}{R}>0\}.$$
Then the field equations on $\mathcal{I}=\{r <R\}$ is redused to the Tolman-Oppenheimer-Volkhoff equation
\begin{equation}
\frac{dm}{dr}=4\pi r^2\rho,\quad
\frac{dP}{dr}=-(\rho+P)\frac{m+4\pi r^3\rho}{r^2(1-\frac{2m}{r})}. \label{TOV}
\end{equation}
Therefore the application of the Cauchy-Kovalevskaja theorem means the shooting of the solution $(m(r), P(r))$ of \eqref{TOV} to the left from $r=R$ with the initial data 
$(m, P)=(M, 0)$. Actually the solution exists locally on $[R-\delta, R]$ with $\delta \ll 1$. The problem is : Can it be prolonged to $r=0$ and hit a regular values at the center? 

In order to fix the idea, we consider the domain $\mathsf{D}$ of the equation \eqref{TOV} as 
\begin{equation}
\mathsf{D}=\{ 0<r, 0<P, 0< m+4\pi r^3\rho , 0<1-\frac{2m}{r} \}.
\end{equation}
(We are assuming that $\rho=\rho(P) >0$ for $P>0$.) Let the left maximal existence of existence of the solution in this specified domain $\mathsf{D}$ to the initial data $(m, P)=(M, 0)$ be $]r_-, R]$. Note that $dP/dr <0$ on $\mathsf{D}$. The following cases might occur:

Either {\bf Case-(0)}: $r_->0$, or {\bf Case-(1)}: $r_-=0$.

If {\bf Case-(0)}, then either
{\bf Case-(00)}: $P \nearrow +\infty$ as $r \searrow r_-$,
or
{\bf Case-(01)}: $P\nearrow \exists  P_-(<\infty)$ as $ r\searrow r_-$.

If {\bf Case-(01)}, it can be shown that $dP/dr \rightarrow 0$ as
$r \rightarrow r_-+0$. 

If {\bf Case-(1)} , either
{\bf Case-(10)}: $P \nearrow +\infty$ as $ r\searrow 0$,
or
{\bf Case-(11)}: $P \nearrow \exists P_{\mathsf{O}} (<\infty)$ as $r \searrow 0$.

Suppose {\bf Case-(11)}. Then 
$$m(r)=M-\int_r^R4\pi\rho(r')(r')^2dr'
\rightarrow \exists m_{\mathsf{O}}
$$
as $r \rightarrow +0$. Since $m(r) <r/2 \rightarrow 0$, we see $m_{\mathsf
{O}}\leq 0$. Since $m(r)+4\pi r^3\rho(r)>0$, we have $m_{\mathsf{O}}\geq 0$.
Therefore $m_{\mathsf{O}}=0$ and
$$m(r)=\int_0^r4\pi\rho(r')(r')^2dr'. $$
Thus if and only if the {\bf Case-(11)}, the solution is that of \eqref{TOV} which is regular at the center. In other words, the Roos' strategy is successful if and only if {\bf Case-(11)}. 

Let us look at the situation in the opposite direction. If we shoot the solution of \eqref{TOV} with the initial data $m=0, P=P_{\mathsf{O}}>0$ to the right from $r=+0$,
the solution $(m(r), P(r))$ may hit $P=0$ at finite $r=R$; when we can verify that 
$M=m(R)$ exists and satisfy $\displaystyle 1-\frac{2M}{R}>0$. The set $\mathcal{O}$ of all such $P_{\mathsf{O}}$ is an open set of $]0,+\infty[$. Of course $\mathcal{O}$ can be empty, which depends on the equation of state $\rho=\rho(P)$. The connected components $\mathcal{O}_j, j=1,2,\cdots,$ are open intervals, and the set of $(R, m(R))$ for 
$P_{\mathsf{O}} \in \mathcal{O}_j$ is a curve $\mathcal{C}_j$ in the set of admissible data $\mathcal{A}$. {\bf Case-(11)} is the case and the Roos' strategy turns out to be successful if and only if $\exists j: (R,M)\in \mathcal{C}_j$. Since the 2-dimensional measure of $\mathcal{C}_j$ is zero, we can say that 
$(R,M) \notin\bigcup_j\mathcal{C}_j$ a. e. In other words {\bf  the strategy of W. Roos is almost everywhere
unsuccessful}. In order to make the Roos' strategy give success we must chose a combination of $M$ and $R$ which fits a very tightly restricted condition, and it seems to be desperate to give an explicit expression of the condition accordingly to the given equation of states.

\section{Ellipsoid as a candidate for the boundary surface}

In the preceding section we examined the strategy of W. Roos for the particular spherically symmetric case. Although it is the case very scarcely, the strategy is successful by taking a sphere $\{r=R\}$ as the matching boundary $\mathcal{S}$. However if the angular velocity $\Omega$ is not equal to $0$, maybe we should take other figures than spheres as $\mathcal{S}$. Inspired by the analogy with the ellipsoidal figures of rotating liquid (that is, incompressible ) stars in the non-relativistic theory (see e.g. \cite{Chandra}), one may try to take ellipsoids as candidates of $\mathcal{S}$ for {\bf (Q2)}. In fact, P. Collas, \cite[p.68]{Collas}, says

\begin{quote}
Hernandez outlined a method for constructing exact interior solution which served as sources for the Kerr metric. The guessed metric matches the Kerr metric on a suitable surface [Foornote (5)] and, in the limit of no rotation, goes into the interior Schwarzschild solution. 

---------

{\it  Footnote (5): A. Krasinski: Institute of Astronomy, Polish Academy of Science preprint No.63, Warsaw (May 1976), has shown that the surface of a source of the Kerr metric should be given by $r=$constant  in Boyer-Lindquist co-ordinates.}
\end{quote}

According to \cite{BoyerL} the Kerr metric is described by the authors R. H. Boyer and R. W. Lindquist as
$$
ds^2= dr^2+2a\sin^2\vartheta dr d\phi +
(r^2+a^2)\sin^2\vartheta d\phi^2 + $$
$$ +\Sigma d\vartheta^2-dt^2+(2Mr/\Sigma)
(dr+a\sin^2\vartheta d\phi+dt)^2, \eqno{\cite[(2.7)]{BoyerL}} $$
with
$$\Sigma=r^2+a^2\cos^2\vartheta. \eqno{\cite[(2.8)]{BoyerL}} $$
But $r$ is a function of $x,y,z$ of the standard co-ordinate system $(t, x, y, z) $ such that
$$\frac{x^2+y^2}{r^2+a^2}+\frac{z^2}{r^2}=1. \eqno{\cite[(2.5)]{BoyerL}} $$
Therefore the surface $r=$ constant referred in the above \cite[Footnote (5)]{Collas} are confocal ellipsoids. (\cite[p.269L, the last line]{BoyerL}.) In other words, \cite{Collas} says that ellipsoids are suitable figures for $\mathcal{S}$ to solve the problem {\bf (Q2)}. \\

However it seems doubtful that an exact ellipsoid is suitable for $\mathcal{S}$, since it seems doubtful that the vacuum boundary $\{r=r_+(\zeta)\}$ of the interior solution constructed in \cite{asEE} would be an exact ellipsoid if $\Omega \not=0$. Let us explain the reason of this doubt in this section.\\

We have
\begin{equation}
r_+(\zeta)=\mathsf{a}\Big[\Xi_1\Big(\zeta,\frac{1}{\gamma-1}, \mathsf{b}\Big)+
O\Big(\frac{u_{\mathsf{O}}^2}{\mathsf{c}^2}\Big)\Big],
\end{equation}
where
$$
\mathsf{a}=\sqrt{\frac{\mathsf{A}\gamma}{4\pi\mathsf{G}(\gamma-1)}}
\rho_{\mathsf{O}}^{-\frac{2-\gamma}{2}},
\qquad
\mathsf{b}=\frac{\Omega^2}{4\pi\mathsf{G}\rho_{\mathsf{O}}}
$$
and
$\displaystyle \xi=\Xi_1\Big(\zeta, \frac{1}{\gamma-1}, \mathsf{b}\Big)$ is the surface curve of the distorted Lane-Emden function
$\Theta\Big(\xi, \zeta, \frac{1}{\gamma-1},\mathsf{b}\Big)$,
while
$$u_{\mathsf{N}}(r,\zeta)=u_{\mathsf{O}}\Theta\Big(\frac{r}{\mathsf{a}}, \zeta, 
\frac{1}{\gamma-1},\mathsf{b}\Big).$$
On the distorted Lane-Emden function, as shown in \cite[\S 7]{JJTM1} and \cite[\S 5]{JJTM2}, we have the approximation
$$\Xi_1(\zeta)=\xi_1+\frac{\xi_1^2}{\mu_1}\mathfrak{h}(\xi_1, \zeta)\mathsf{b}+
O(\mathsf{b}^{\frac{1}{\gamma-1}\wedge 2}),$$
where $\xi_1=\xi_1(\frac{1}{\gamma-1}), \mu_1=\mu_1(\frac{1}{\gamma-1})$ are positive numbers,
$$\mathfrak{h}(\xi,\zeta)=h_0(\xi)+A_2\psi_2(\xi)P_2(\zeta), $$
$A_2<0, \psi_2(\xi_1)>0$, and
$$P_2(\zeta)=\frac{1}{2}(3\zeta^2-1). $$
Thus $r=r_+(\zeta)$ is approximated by
$$r=\mathsf{a}\Big[c_0+(c_1-c_2\zeta^2)\mathsf{b}\Big]$$
with $c_2>0$, which does not give an ellipsoid if $\mathsf{b}\not=0$, since an ellipsoid should be given by a function of the form
$$r =\frac{a_0}{\sqrt{1+a_1\zeta^2}}.$$ Of course this is not a rigorous proof of the claim that $r=r_+(\zeta)$ does not give an ellipsoid, since the detailed structure of the remainder terms $O(u_{\mathsf{O}}^2/\mathsf{c}^2)$ and 
$O(\mathsf{b}^{\frac{1}{\gamma-1}\wedge 2})$ are not clearly analyzed. But it is strongly plausible that $\xi=\Xi_1(\zeta)$ does not give an ellipsoid if $\mathsf{b}\not=0$. Actually the following ``no-go
theorem'' has been known for more than 120 years:\\

{\bf Theorem of Hamy-Pizzetti: An ellipsoidal stratification} ({\it that is, the situation that all level surfaces 
$\{u =\mbox{Const.}\}$ are ellipsoids }) {\bf is impossible for heterogeneous} ({\it
that is, with non-constant $\rho$} ), {\bf rotating symmetric figures of equilibrium} ({\it that is, axisymmetric stationary solutions of the Euler-Poisson equations with 
$\Omega\not=0$}).\\

See \cite[Sec. 3.2]{Moritz}, \cite{Pizzetti}, \cite{Wavre}. By this no-go theorem it is impossible that all level surfaces $\{\rho =\mbox{Const.}\}$ are ellipsoids, since
$\partial\Theta/\partial \xi <0$ for $0<\xi <\Xi_1(\zeta)$ so that $\partial\rho/\partial r <0$ for $0<r<r_+(\zeta)$ in the non-relativistic problem for which $\mathsf{c}=\infty$,
$ r_+(\zeta)=\mathsf{a}\Xi_1(\zeta)$, and $$\rho=\rho_{\mathsf{N}}=
\Big(\frac{\gamma-1}{\mathsf{A}\gamma}\Big)^{\frac{1}{\gamma-1}}
(u_{\mathsf{N}}\vee 0)^{\frac{1}{\gamma-1}}
=\rho_{\mathsf{O}}(\Theta \vee 0)^{\frac{1}{\gamma-1}},$$
although this no-go theorem does not claim that it is impossible that the individual level surface $\{ u (=u_{\mathsf{N}} )=0\} =\{ r=r_+(\zeta) (=\mathsf{a}\Xi_1(\zeta) ) \}$ is an ellipsoid.

\end{document}